\documentclass{amsart2000}
\usepackage{amscd2000,amsmath2000,amssymb,subfigure,hyperref}
\usepackage[arrow,matrix,graph,frame,poly,arc,tips]{xy}

\theoremstyle{plain}
\newtheorem{thm}[subsection]{Theorem}
\newtheorem{lem}[subsection]{Lemma}

\newtheorem{cor}[subsection]{Corollary}

\theoremstyle{definition}
\newtheorem{remark}[subsection]{Remark}
\newtheorem{definition}[subsection]{Definition}
\newtheorem{example}[subsection]{Example}
\newtheorem{conj}[subsection]{Conjecture}

\numberwithin{equation}{section}

\newcommand{\DS}{\displaystyle }
\newcommand{\abs}[1]{\lvert #1 \rvert}

\newcommand{\A}{{\mathcal A}}
\newcommand{\B}{{\mathcal B}}

\newcommand{\Z}{{\mathbb Z}}

\newcommand{\C}{{\mathbb C}}
\newcommand{\G}{\mathsf{G}}
\renewcommand{\k}{\Bbbk}

\newcommand{\m}{\mathfrak{m}}

\DeclareMathOperator{\Tor}{Tor}
\DeclareMathOperator{\Ext}{Ext}

\DeclareMathOperator{\Hilb}{Hilb}
\DeclareMathOperator{\rank}{rank}

\DeclareMathOperator{\gr}{gr}

\DeclareMathOperator{\hor}{hor}
\DeclareMathOperator{\ver}{vert}
\DeclareMathOperator{\tot}{tot}
\DeclareMathOperator{\im}{im}
\DeclareMathOperator{\coker}{coker}
\DeclareMathOperator{\codim}{codim}
\DeclareMathOperator{\spn}{span}
\newcommand{\Kosz}{\mathcal K}

\begin{document}

\title[Lower central series and free resolutions of arrangements]%
{Lower central series and free resolutions of\\ hyperplane arrangements}

\author[Henry K. Schenck]{Henry K. Schenck$^1$}
\address{Department of Mathematics,
Harvard University, Cambridge, MA 02138}
\curraddr{Department of Mathematics,
Texas A\&M University, College Station, TX 77843}
\email{\href{mailto:schenck@math.tamu.edu}{schenck@math.tamu.edu}}
\urladdr{\href{http://www.math.tamu.edu/~schenck/}%
{http://www.math.tamu.edu/\~{}schenck}}

\thanks{$^1$Partially supported by an NSF postdoctoral research fellowship}

\author[Alexander I. Suciu]{Alexander~I.~Suciu$^2$}
\address{Department of Mathematics,
Northeastern University,
Boston, MA 02115}
\email{\href{mailto:alexsuciu@neu.edu}{alexsuciu@neu.edu}}
\urladdr{\href{http://www.math.neu.edu/~suciu/}%
{http://www.math.neu.edu/\~{}suciu}}

\thanks{$^2$Partially supported by NSF grant DMS-0105342}




\begin{abstract}
If $M$ is the complement of a hyperplane arrangement, and $A=H^*(M,\k)$ is the
cohomology ring of $M$ over a field of characteristic $0$, then the ranks,
$\phi_k$, of the lower central series quotients of $\pi_1(M)$ can be computed
from the Betti numbers, $b_{ii}=\dim \Tor^A_i(\k,\k)_i$,
of the linear strand in a minimal free resolution of $\k$ over $A$.
We use the Cartan-Eilenberg change of rings spectral sequence to relate
these numbers to the graded Betti numbers, $b'_{ij}=\dim \Tor^E_i(A,\k)_{j}$,
of a minimal resolution of $A$ over the exterior algebra $E$.

 From this analysis, we recover a formula of Falk for $\phi_3$,
and obtain a new formula for $\phi_4$. The exact sequence of
low degree terms in the spectral sequence allows us to answer
a question of Falk on graphic arrangements, and also shows that
for these arrangements, the algebra $A$ is Koszul if and only if
the arrangement is supersolvable.

We also give combinatorial lower bounds on the Betti
numbers, $b'_{i,i+1}$, of the linear strand  of the free resolution of
$A$ over $E$; if the lower  bound is attained for $i=2$, then it is
attained for all $i \ge 2$. For such arrangements, we compute the entire
linear strand of the resolution, and we prove that all components of
the first resonance variety of $A$ are local.  For graphic arrangements
(which do not attain the lower bound, unless they have no braid
sub-arrangements), we show that  $b'_{i,i+1}$ is determined by
the number of triangles and $K_4$ subgraphs in the graph.
\end{abstract}

\maketitle


\section{Introduction}
\label{sec:intro}

\subsection{LCS formulas}
\label{subsec:introLCS}
The lower central series (LCS) of a finitely-generated group $G$ is
a chain of normal subgroups, $G=G_1 \ge G_2 \ge G_3\ge\cdots$,
defined inductively by $G_k = [G_{k-1},G]$.
Fix a field $\k$ of characteristic $0$.
Associated to $G$ is a graded Lie algebra,
$\gr(G)\otimes \k := \bigoplus_{k=1}^{\infty} G_k/G_{k+1} \otimes \k$,
with Lie bracket induced by the commutator map.
Let $\phi_k=\phi_k(G)$ denote the rank of the $k$-th quotient.

One case in which the numbers $\phi_k$ have aroused great
interest occurs when $M$ is the complement of a complex
hyperplane arrangement and $G = \pi_1(M)$ is its fundamental group.
In certain situations, the celebrated Lower Central Series formula holds:
\begin{equation}
\label{lcsformula}
\prod_{k=1}^{\infty}(1-t^k)^{\phi_k}=P(M,-t),
\end{equation}
where $P(M,t)=\sum b_it^i $ is the Poincar\'e polynomial of the complement.
This formula was first proved by Kohno \cite{K} for braid arrangements, and
then generalized by Falk and Randell \cite{FR85} to supersolvable arrangements.
In \cite{SY}, Shelton and Yuzvinsky gave a beautiful interpretation of
the LCS formula in terms of Koszul duality. The LCS formula  holds for
rational $K(\pi,1)$ spaces (see Papadima and Yuzvinsky \cite{PY}), and
an analogue holds for hypersolvable arrangements (see Jambu and 
Papadima \cite{JP}).

All these LCS formulas ultimately hinge on the assumption that the
cohomology algebra (or its quadratic closure) is Koszul.
The goal of this paper is to explore formulas of type \eqref{lcsformula}
for arrangements where such assumptions fail.  In complete generality,
there is no good substitute for the polynomial $P(X,-t)$ on the
right side, as illustrated by the non-Fano plane (see Peeva \cite{Pe}).
But, as conjectured in \cite{Su}, and as verified here in certain
situations (at least in low degrees), there do exist large classes of
arrangements for which a modified LCS formula holds, with the
Poincar\'e polynomial replaced by another (combinatorially
determined) polynomial.

\subsection{LCS ranks and Betti numbers of the linear strand}
\label{subsec:ext}

Our starting point is a formula expressing the
LCS ranks of an arrangement group in terms of the
Betti numbers,
$b_{ii}= \dim_{\k} \Ext^i_A(\k,\k)_i=\dim_{\k} \Tor_i^A(\k,\k)_i$,
of the linear strand of the free resolution of the residue field $\k$,
viewed as a module over the Orlik-Solomon algebra, $A=H^*(M,\k)$:
\begin{equation}
\label{eq:syintro}
\prod_{k=1}^{\infty} (1-t^k)^{-\phi_k} =
\sum_{i=0}^{\infty} b_{ii}t^i.
\end{equation}

A bit of history:
In \cite{K}, Kohno showed that the left hand quantity is equal to the 
Hilbert series
of the universal enveloping algebra, $U=U(\mathfrak{g})$, of the 
holonomy Lie algebra
$\mathfrak{g}$ of the complement $M$.
(This follows from the Poincar\'e-Birkhoff-Witt theorem, together with the
isomorphism $\gr(G)\otimes\k\cong\mathfrak{g}$, which is a consequence of
the formality of $M$, in the sense of Sullivan.)    In \cite{SY},
Shelton and Yuzvinsky proved that $U=\overline{A}^{!}$, the Koszul dual of the
quadratic closure of the Orlik-Solomon algebra. Finally, results of
Priddy \cite{Pr} and L\"ofwall~\cite{L} relate the Koszul dual of a
quadratic algebra to the linear strand in the Yoneda Ext-algebra:
\begin{equation}
\label{eq:shriek}
\overline{A}^{!} \cong \bigoplus_i\Ext^i_{\overline{A}}(\k,\k)_i.
\end{equation}
Since obviously $\Ext^i_{\overline{A}}(\k,\k)_i=\Ext^i_A(\k,\k)_i$,
formula \eqref{eq:syintro} follows at once (compare \cite[Theorem~2.6]{Pe}).

If $A$ is a Koszul algebra
(i.e., $\Ext^i_A(\k,\k)_j=0$, for $i\ne j$), then $A=\overline{A}$,
and the Koszul duality formula, $\Hilb(A^{!},t)\cdot \Hilb(A,-t)=1$,
yields the LCS formula \eqref{lcsformula}. This computation applies 
to supersolvable
arrangements, which do have Koszul OS-algebras, cf.~\cite{SY}.
The more general formula \eqref{eq:syintro} was first exploited by Peeva in
\cite{Pe} to obtain bounds on the LCS ranks $\phi_k$.

\subsection{The Orlik-Solomon algebra and free resolutions}
\label{subsec:introOS}
A minimal free resolution of a graded module $N$ over a graded $\k$-algebra
$R$ is simply a graded exact sequence
\begin{equation}
\label{eq:freeres}
P_{\bullet}\colon \ \cdots \longrightarrow \bigoplus_{j}  R^{b_{2j}}(-j)
\longrightarrow
\bigoplus_{j} R^{b_{1j}}(-j) \longrightarrow \bigoplus_{j} R^{b_{0j}}(-j)
\longrightarrow N\to 0,
\end{equation}
where the exponents are
$b_{ij}(N)=\dim_{\k} \Tor_i^R(N,\k)_j$.
We are especially interested in  two cases:
\begin{itemize}
\item $N=A=E/I$, the Orlik-Solomon algebra of an arrangement
of $n$ hyperplanes, and $R=E$ is the exterior algebra on $n$ generators,
in which case we will write $b'_{ij}=\dim_{\k} \Tor_i^E(A,\k)_j$.

\item $N=\k$ and $R=A$,
in which case we will write $b_{ij}=\dim_{\k} \Tor_i^A(\k,\k)_j$.
\end{itemize}

We study the relationship between these two free resolutions by
means of the change of rings spectral sequence associated to the
composition $E \rightarrow A \twoheadrightarrow \k$:
\begin{equation}
\label{eq:ssintro}
\Tor_i^A\left(\Tor_j^E(A,\k),\k\right)
\Longrightarrow \Tor_{i+j}^E(\k,\k).
\end{equation}
The goal is to translate the right-side of formula \eqref{eq:syintro}
into an expression involving the resolution of $A$ over $E$,
which is a much smaller
(in the sense that the ranks of the free modules appearing in the
resolution are much smaller) object that the resolution of $\k$
over $A$.

A resolution $P_{\bullet}$ as in \eqref{eq:freeres} is {\it linear}
if each free module $F_i$ is generated in a single degree
$\alpha_i$ and $\alpha_{i} = \alpha_{i-1}-1$.
By convention, when $N=R/J$ and $J$
is generated in degree $\ge 2$, the {\it linear strand} is the
complex with terms $F_i=R^{b_{i,i+1}}(-i-1)$, $i \ge 1$.
Most of the results we obtain from the change of rings spectral sequence
stem from the happy fact that the resolution of $\k$ as a module over
$E$ is linear (with $\dim_{\k} \Tor_i^E(\k,\k)_i = \binom{n + i -1}{i}$).
In \cite{EPY}, Eisenbud, Popescu and Yuzvinsky prove that
$A^{\vee}=H_*(M,\k)$, viewed as an $E$-module via the cap product,
has a linear free resolution;
we had hoped  that this would make one of the change of rings
spectral sequences particularly simple,  but it does not.

\subsection{LCS ranks of arrangement groups}
\label{subsec:lcsarr}
Let $a_i$ denote the number of minimal
generators of degree $i$ in the Orlik-Solomon ideal $I$.
It is elementary to see that $\phi_1=b_1$ and $\phi_2=a_2$.
In \cite{f}, Falk gave a formula for
$\phi_3$ of an arbitrary arrangement, which can be
interpreted as $\phi_3=b'_{23}$,  but no general formula for
$\phi_4$ was known (besides the one implicit in \eqref{eq:syintro}).
Using the above spectral sequence, we obtain a formula for $\phi_4$
which depends solely on the resolution of $A$ over $E$:
\begin{equation}
\label{eq:phi4intro}
\phi_4=\binom{a_2}{2}+b'_{34}-\delta_4,
\end{equation}
where $\delta_4$ denotes the number of minimal quadratic
syzygies on the degree~$2$ generators of $I$ which are Koszul syzygies.
For many arrangements, this expression yields an explicit combinatorial formula
for $\phi_4$.

\subsection{Minimal linear strand arrangements}
\label{subsec:loclcs}
For arrangements with
\begin{equation}
\label{eq:alocal}
b'_{23}=2\sum _{X\in L_2(\A)} \binom{\mu(X)+1}{3},
\end{equation}
we have
\begin{equation}
\label{eq:nrphi4}
\phi_4=\sum _{X\in L_2(\A)} \frac{\mu(X)^2 (\mu(X)^2-1)}{4}.
\end{equation}

Arrangements for which \eqref{eq:alocal} holds have $b'_{i,i+1}$ as small as
possible, for all $i\ge 2$ (given the combinatorics); in other words,
the linear strand in the free resolution of $A$ over $E$ is minimal.
In view of this, we call such arrangements {\em minimal linear strand (MLS)}
arrangements.   Examples include Kohno's arrangements
$X_2$ and $X_3$,  the non-Pappus arrangement,
and graphic arrangements with  no $K_4$ subgraphs.
For MLS arrangements in $\C^3$, we compute the whole Betti diagram
of the minimal free resolution of $A$ over $E$, as  well as the
differentials in the linear strand.

As a generalization of formula \eqref{eq:nrphi4}, we conjecture that,
for all MLS arrangements
and all $k\ge 2$,
\begin{equation}
\label{eq:mobiusintro}
\phi_k=\tfrac{1}{k} \sum_{d\mid k}\sum_{X\in L_2(\A)}
\mu(d) \mu(X)^{\frac{k}{d}},
\end{equation}
where $\mu$ stands for both the
M\"obius function of the integers,
and the M\"obius function of the
intersection lattice of the arrangement.
Equivalently:
\begin{equation}
\label{eq:mob}
\prod_{k=1}^{\infty}(1-t^k)^{\phi_k}=
(1-t)^{b_1-b_2}\prod_{X\in L_2(\A)} (1- \mu(X)\, t).
\end{equation}
The proof of the conjecture hinges upon computing the image of a
$d_2$ map in the change of rings spectral sequence.

\subsection{Graphic arrangements}
\label{subsec:graph}
As another application our methods, we study the
free resolutions and LCS quotients of a graphic arrangement.
Let $\G$ be a (simple) graph on $\ell$ vertices, with edge-set $\mathsf{E}$,
and let $\A_{\G}=\{z_i-z_j=0\mid (i,j)\in \mathsf{E} \}$
be the corresponding arrangement in $\C^{\ell}$.
Using a well-known result of Stanley, we show that the OS-algebra
$A$ is Koszul if and only if the arrangement $\A_{\G}$ is supersolvable.

As for the LCS ranks of the fundamental group of the complement
of $\A_{\G}$, we conjecture that they are given by
\begin{equation}
\label{eq:kap}
\phi_k=\tfrac{1}{k} \sum_{d\mid k}
\sum_{j=1}^{k}\sum_{s=j}^{k} (-1)^{s-j} \binom{s}{j} \kappa_{s}
   \mu(d) j^{\frac{k}{d}}
\end{equation}
where $\kappa_s$ is the number of complete subgraphs on $s+1$ vertices;
equivalently:
\begin{equation}
\label{eq:kappaintro}
\prod_{k=1}^{\infty} \left(1-t^k\right)^{\phi_k}=
\prod_{j=1}^{\ell-1}
\left(1-jt\right)^{\DS{\sum_{s=j}^{\ell-1}(-1)^{s-j} \tbinom{s}{j} \kappa_{s}
}}
\end{equation}
We verify this conjecture for $k\le 3$. In particular, we find
$\phi_3=2(\kappa_2+\kappa_3)$, thereby answering a question of
Falk \cite{Fa00}.  For $k=4$, we verify the inequality
$\phi_4\ge 3(\kappa_2+3\kappa_3+2\kappa_4)$, and show
that equality holds if $\kappa_3=0$.

\subsection*{Acknowledgment}
The computations for this work were done primarily
with the algebraic geometry and commutative algebra system
{\sl Macaulay~2}, by Grayson and Stillman~\cite{GS}.
Additional computations  were done with {\sl Mathematica~4.0},
from Wolfram Research.

\section{Resolution of the OS-algebra over the exterior algebra}
\label{sec:resos}

We start by analyzing the minimal free resolution of the
Orlik-Solomon algebra over the exterior algebra.  For
``minimal linear strand" arrangements, the graded Betti
numbers of this resolution can be computed explicitly from
the M\"obius  function of the intersection lattice.

\subsection{Orlik-Solomon algebra}
\label{subsec:osalg}
Let $\A$ be a (central) arrangement of complex hyperplanes
in $\C^{\ell}$, with complement $M(\A)=\C^{\ell}\setminus
\bigcup_{H\in \A} H$.
Let $L(\A)=\{\bigcap_{H\in \B}\, H \mid \B\subseteq \A\}$ be the
{\em intersection lattice} of all flats of the arrangement,
with order given by reverse inclusion, and rank given by codimension.
The cohomology ring of $M(\A)$ with coefficients in a field $\k$
is isomorphic to a certain $\k$-algebra, $A=A(\A)$, which can be described
solely in terms of $L(\A)$, as follows (see the book by Orlik-Terao
\cite{ot} as a general reference).

Let $E$ be the exterior algebra on generators $e_1,\dots,e_n$
in degree $1$, corresponding to the hyperplanes.
This a differential graded algebra over $\k$,
with  grading
\[
E_r=\spn_{\k}\, \{ e_J=e_{i_1}\wedge \cdots \wedge e_{i_r} \mid
J=\{i_1,\dots ,i_r\}\subseteq \{1,\dots,n\}\},
\]
and differential $\partial\colon E_i\to E_{i-1}$ given by
$\partial(1)=0$, $\partial(e_i)=1$, together with the
graded Leibnitz rule. The {\em Orlik-Solomon algebra}, $A$, is
the quotient of  $E$ by the (homogeneous) ideal
\[
I=\Big(\partial e_{\B}\, \big|\: \codim \bigcap_{H\in \B} H < |\B| \big.\Big).
\]

The algebra $A=E/I$ inherits a grading from $E$ in the obvious fashion.
We denote the Betti
numbers of the arrangement by $b_i=\dim A_i$.  As is well-known,
\begin{equation}
\label{eq:betti}
b_i=\sum_{X\in L_i(\A)} (-1)^i \mu(X),
\end{equation}
where $\mu\colon L(\A)\to \Z$ is the M\"{o}bius function of the
lattice.

Let $a_j$ be the number of minimal generators of $I$ having degree
exactly $j$, and let $I[j] \subseteq I$ be the ideal of $E$ generated
by those $a_j$ elements. We have:
\begin{equation}
\label{eq:aj}
a_j = \dim_{\k} (I[j])_j = \dim_{\k} \Tor_1^E(A,\k)_j.
\end{equation}
  From the exact sequence
$0\to I \to E \to A \to 0$, we get:
$a_j+ b_j\le\binom{n}{j}$. Furthermore, since
$I$ is generated in degree $\ge 2$, we have $a_1=0$ and
\begin{equation}
\label{eq:a2}
a_2=\binom{n}{2}-b_2=\sum_{X\in L_2(\A)} \binom{\mu(X)}{2}.
\end{equation}
If $A$ is quadratic (i.e., $I=I[2]$), then
$a_j=0$ for $j>2$.  No explicit combinatorial formula
for the numbers $a_j$ ($j>2$) is known in general.

\subsection{Free resolution of the OS-algebra}
\label{subsec:osres}
By general results on modules over the exterior algebra,
the Orlik-Solomon algebra $A$ admits a resolution over $E$
by finitely-generated free modules,
\begin{equation}
\label{eq:ares}
\cdots \longrightarrow F_2 \longrightarrow F_1\longrightarrow E
\longrightarrow A\longrightarrow 0.
\end{equation}
Furthermore, there exist {\em minimal} free resolutions,
in the sense that $\dim_{\k} F_i =\dim_{\k} \Tor_i^E(A,\k)$,
for all $i\ge 1$. Such resolutions are unique up to chain-equivalence.
Our goal in this section is to obtain a better understanding of the
graded Betti numbers
\begin{equation}
\label{eq:bettiprime}
b'_{ij}=\dim_{\k} \Tor_i^E(A,\k)_j \, .
\end{equation}
We start with some simple observations.

\begin{lem}
\label{lem:torvanish}
Let $\A$ be an essential, central arrangement of rank $\ell$.
If $j \ge i+\ell$, then
\begin{equation*}
\Tor_i^E(A,\k)_j = 0.
\end{equation*}
\end{lem}

\begin{proof}
As in the proof of Theorem~1.1 of \cite{EPY}, to study the free $E$-resolution
of $A$, we may pass to a decone of $\A$.
Let $A'=E'/I'$ be the  OS-algebra of the decone.
Then $A'_j=0$, for $j \ge \ell$.  The free $E'$-resolution of
$\k$ is  linear, so after tensoring this resolution with $A'$,
we obtain a complex of modules with
$i$-th term $\bigoplus A'(-i)$,
which vanishes in degrees $i + \ell$ and higher.
\end{proof}

\begin{lem}
\label{lem:bb}
Let $\A$ be an arrangement of $n$ hyperplanes in $\C^{\ell}$.  The 
Betti numbers
of the minimal free resolution of $A$ over $E$ are related
to the Betti numbers of $A$, as follows:
\begin{equation}
\label{eq:bprime}
\Bigg(\sum_{i=0}^{\infty}\sum_{j=0}^{\ell-1} (-1)^i
b'_{i,i+j}t^{i+j} \Bigg) \cdot
(1+t)^n=\sum_{i=0}^{\ell} b_it^i
\end{equation}
\end{lem}

\begin{proof}
  From \eqref{eq:ares}, we get
$\big(\sum_{i}(-1)^i \Hilb(F_i,t)\big)\cdot P(E,t)=P(A,t)$.
\end{proof}

\subsection{Local syzygies}
\label{subsec:locsyz}
We now analyze in more detail a certain type of linear syzygies
in the resolution of $A$ over $E$.  These ``local" syzygies
come from flats in $L'_2(\A):=\{X\in L_2(\A)\mid \mu(X)>1\}$.
We start with the simplest situation.

\begin{example}[Pencil of $3$ lines]
\label{ex:pen3}
Let $\A=\{H_0,H_1,H_2\}$ be an arrangement of $3$ lines through the
origin of $\C^2$.  The Orlik-Solomon ideal is
generated by $\partial e_{012}=(e_1-e_2)\wedge (e_0-e_2)$.
It is readily seen that the minimal free resolution of
$A$ over $E$ is:
{\Small
\begin{equation*}
\begin{split}
0 \leftarrow A \longleftarrow E
\xleftarrow{\begin{pmatrix}\partial e_{012}\end{pmatrix}}
F_1
\xleftarrow{\begin{pmatrix}
e_1-e_2 & e_0-e_2
\end{pmatrix}}
F_2
\xleftarrow{\begin{pmatrix}
e_1-e_2 & e_0-e_2 & 0\\0& e_1-e_2&e_0-e_2
\end{pmatrix}}
F_3
\longleftarrow\\[4pt]
\cdots \longleftarrow
F_i
\xleftarrow{\begin{pmatrix}
\xymatrixrowsep{-1pt}
\xymatrixcolsep{2pt}
\xymatrix{
e_1-e_2 & e_0-e_2 & 0  & \cdots &0 & 0\\
0 & e_1-e_2\ar@{..}[ddrrr]
& e_0-e_2\ar@{..}[ddrrr]   & \cdots &0 & 0\\
& &&&&\\
0& 0 & 0 &\cdots &e_1-e_2&e_0-e_2
}
\end{pmatrix}}
F_{i+1} \longleftarrow \cdots
\end{split}
\end{equation*}
}

\noindent
with $F_i=E^i(-i)$.
Thus, $b'_{i,i+1}=i$, for $i\ge 1$, and $b'_{i,i+r}=0$, for $r>1$.
\end{example}

\begin{example}[Pencils of $n+1$ hyperplanes]
\label{ex:pencils}
More generally, let $\A=\{H_0,\dots,H_n\}$ be a central arrangement of
$n+1\ge 3$ hyperplanes whose common intersection is of codimension two.
The Orlik-Solomon ideal is generated by all elements
$\partial e_{ijk}=(e_j-e_i)\wedge (e_k-e_i)$ with $0\le i<j<k\le n$.
In fact, $I=I[2]$, with minimal set of generators
$(e_j-e_0)\wedge (e_k-e_0)$, for $1\le j<k\le n$.

The minimal resolution of $A$ over $E$ has a similar form to the one above.
The graded Betti numbers are given by:
\begin{equation}
\label{bettipencil}
b'_{i,i+1}=i\,\binom{n+i-1}{i+1}\, , \quad \text{for all $i\ge 1$},
\end{equation}
and $b'_{i,i+r}=0$, for $r>1$.  To see this, note that,
after a suitable linear change of variables, the ideal $I$
becomes identified with the  monomial ideal
$\mathfrak{m}^2$, where $\mathfrak{m}=\left(e_1,\dots,e_n\right)$.
Formula \eqref{bettipencil}  then follows
from general results of \cite{AAH}, or directly, from the exact sequence
\begin{equation*}
0 \longrightarrow \Tor^{E'}_{i+1}(\k,\k) \longrightarrow
\Tor^{E'}_i(\m/\m^2(-1),\k)
\longrightarrow \Tor^{E'}_i(E'/\m^2,\k) \longrightarrow 0,
\end{equation*}
where $E'=E/(e_0)$, together with the isomorphism $\m/\m^2 \cong \k^n$.
\end{example}

Going back to an arbitrary arrangement, we observe that
for each flat $X\in L_2(\A)$ with $\mu(X)\ge 2$ there are linear
syzygies of the type discussed in the examples above,
which we will call {\em local} linear syzygies.
In fact, these syzygies persist in the free
resolution of $A$ over $E$, as the next lemma shows.

\begin{lem}
\label{lem:bprimebound}
Let $\A$ be an arrangement. Then:
\begin{equation}
\label{eq:bpineq}
b'_{i,i+1}\ge i\sum _{X\in L_2(\A)} \binom{\mu(X)+i-1}{i+1}.
\end{equation}
Moreover, if equality
holds for $i=2$, then it holds for all $i\ge 2$.
\end{lem}

\begin{proof}
Sets of local linear syzygies corresponding to different elements
of $L_2(\A)$ are linearly independent. Indeed, they are supported
only on generators coming from the same element of $L_2(\A)$.
\end{proof}

\subsection{MLS arrangements}
\label{subsec:locarr}

In view of the preceding lemma, we single out a class of arrangements
for which the linear strand of a minimal free resolution
of $A$ over $E$ is completely determined by the
M\"obius function of $L_2(\A)$.

\begin{definition}
An arrangement is called {\em minimal linear strand (MLS)} if
\begin{equation}
\label{a3cond}
b'_{23}=2\sum _{X\in L_2(\A)} \binom{\mu(X)+1}{3}.
\end{equation}
\end{definition}

In \cite{EPY}, Eisenbud, Popescu, and Yuzvinsky remark that
it is an interesting problem to determine explicitly a minimal
free resolution of the OS-algebra over the exterior algebra.
In the case of MLS arrangements in $\C^3$, Lemmas \ref{lem:bb} and
\ref{lem:bprimebound} permit us to compute the whole Betti diagram
of the minimal resolution of $A$ over $E$, as  well as the
differentials in the linear strand.

\begin{thm}
\label{thm:minfreeres}
Let $\A$ be an MLS arrangement of $n$ hyperplanes in $\C^3$.
Then, the graded Betti numbers of a minimal free resolution of the OS-algebra
over the exterior algebra are given by:
\begin{align}
\label{eq:bp1}
b'_{i,i+1}&= i\sum _{X\in L_2(\A)} \binom{\mu(X)+i-1}{i+1}\\
\label{eq:bp2}
b'_{i,i+2}&= b'_{i+1,i+2}
+\binom{i+1}{2}\binom{n+i-2}{i+2}
-\binom{n+i-2}{i} \sum _{X\in L_2(\A)} \binom{\mu(X)}{2}
\end{align}

\end{thm}

\begin{proof}
The formula for $b'_{i,i+1}$ is given by Lemma~\ref{lem:bprimebound}.
The formula for $b'_{i,i+2}$ follows from \eqref{eq:bprime}, which
in the case $\ell=3$ boils down to:
\begin{equation*}
\label{eq:bprime2}
\qquad
1 -a_2t^2 +\sum_{i=1}^{\infty}
(-1)^i (b'_{i,i+2}-b'_{i+1,i+2})t^{i+2} =
\frac{1+(n-1)t+\big(\binom{n-1}{2}-a_2\big) t^2}{(1+t)^{n-1}}.
\end{equation*}
The differentials in the linear strand are obtained in the obvious
fashion from the differentials in Examples \ref{ex:pen3} and \ref{ex:pencils}.
\end{proof}

\begin{figure}
\subfigure{%
\begin{minipage}[t]{0.3\textwidth}
\setlength{\unitlength}{16pt}
\begin{picture}(5,6)(0,-1.2)
\put(2,2){\oval(5,5)[t]}
\put(2,1.5){\oval(5,5)[b]}
\multiput(-0.52,1.5)(5,0){2}{\line(0,1){0.5}}
\multiput(0,1)(0,2){2}{\line(1,0){4}}
\multiput(1,-0.6)(2,0){2}{\line(0,1){4.6}}
\put(0.34,4){\line(2,-3){3.05}}
\end{picture}
\end{minipage}
}
\setlength{\unitlength}{0.7cm}
\subfigure{%
\begin{minipage}[t]{0.3\textwidth}
\begin{picture}(3,2.5)(0,-0.7)
\put(3,3){\line(1,-1){3}}
\put(3,3){\line(-1,-1){3}}
\put(0,0){\line(1,0){6}}
\multiput(0,0)(6,0){2}{\circle*{0.3}}
\multiput(1.5,1.5)(3,0){2}{\circle*{0.3}}
\multiput(3,3)(0,-3){2}{\circle*{0.3}}
\end{picture}
\end{minipage}
}
\caption{\textsf{The arrangement $\operatorname{X}_3$
and its associated matroid}}
\label{fig:x3}
\end{figure}

\begin{example}
\label{ex:x3}
Let $\A$ be the $\text{X}_3$ arrangement, with defining
polynomial $Q=xyz(x-z)(y+z)(y+2x)$, see Figure~\ref{fig:x3}.   We have
$L'_2=\{X_1,X_2,X_3\}$, with $\mu(X_i)=2$, and so $a_2=3$.
Moreover, $a_3=1$, and so condition \eqref{a3cond}
is satisfied.  From Theorem~\ref{thm:minfreeres}, we compute:
\begin{align*}
b'_{i,i+1}&=3i
\\
b'_{i,i+2}&=\frac{i(i+1)(i^2+5i-2)}{8}
\end{align*}

More generally, for each $n\ge 6$, consider the arrangement
with defining polynomial
$Q=xyz(x-z)(y+z)(y+2x)\cdots (y+(n-4)x)$. We then have
$L'_2=\{X_1,X_2,X_3\}$, with $\mu(X_1)=n-4$, $\mu(X_2)=\mu(X_3)=2$,
and $a_3=n-5$, and so condition \eqref{a3cond}
is satisfied. We leave the computation of the Betti numbers $b'_{ij}$
as an exercise.
\end{example}

\subsection{A criterion for quadraticity}
\label{subsec:quadcrit}

 From the proof of Theorem~\ref{thm:minfreeres}, we see that:
\begin{equation}
\label{eq:bp23}
b'_{23}-a_3=b_3-\binom{b_1}{3}+b_1a_2.
\end{equation}
Applying Lemma~\ref{lem:bprimebound}, we get:
\begin{equation}
\label{eq:quadtest}
a_3\ge \binom{b_1}{3}-b_3-b_1a_2 +2\sum _{X\in L_2(\A)} \binom{\mu(X)+1}{3}.
\end{equation}
This inequality, when coupled with the definition of $a_3$, provides
a necessary
combinatorial criterion for quadraticity  (and hence, Koszulness) of
the OS-algebra.

\begin{thm}
\label{thm:quadcrit}
If $\binom{b_1}{3}-b_3-b_1a_2 +2\sum _{X\in L_2(\A)} \binom{\mu(X)+1}{3}>0$,
then $A$ is not quadratic.
\end{thm}

This criterion shows that  none of the arrangements
in Example~\ref{ex:x3} has a quadratic OS-algebra.

\section{Resolution of the residue field over the OS-algebra}
\label{sec:freeres}

We now turn to the minimal free resolution of the residue field, $\k$,
over the Orlik-Solomon algebra, $A$, and to the corresponding
graded Betti numbers,
$b_{ij}=\dim_{\k}\Tor^A_i(\k,\k)_j$.  A $5$-term exact sequence argument
relates $b_{2j}$ to the number of  minimal generators in degree $j$ of the
Orlik-Solomon ideal, $a_j = \dim_{\k} (I[j])_j$.  This computes the
first three LCS ranks $\phi_k$.

\subsection{Relating various Betti numbers}
\label{subsec:linbetti}
Since the Hilbert series of the residue field is simply $1$,
there are simple numerical constraints on the free resolution of $\k$ 
over $A$.
Suppose $\A$ is a central, essential arrangement in $\C^{\ell}$.  Then:
\begin{equation}
\label{eq:hilb1}
1 = \big(1 - b_{11}t +\sum\limits_{j}b_{2j}t^j
-\sum\limits_{j}b_{3j}t^j +\cdots \big)
(1+ b_{1}t + b_{2}t^2+\cdots + b_{\ell}t^{\ell}).
\end{equation}
This gives us a way to solve for $b_{rr}$, if we know
all $b_{ir}$ with $i<r$:
\begin{center}
$\begin{array}{ccc}
b_{11} &=& b_1\\
b_{22} &=& b_1^2-b_2\\
b_{33} - b_{23} &=& b_1^3-2b_1b_2+b_3\\
b_{44} - b_{34} + b_{24}&=& b_1^4-3b_1^2b_2+b_2^2+2b_1b_3-b_4\\
b_{55} - b_{45} + b_{35}-b_{25}&=&
b_{1}^5 - 4b_{1}^3b_{2} + 3b_{1}b_{2}^2 + 3b_{1}^2b_{3} -
    2b_{2}b_{3} - 2b_{1}b_{4} + b_{5}
\end{array}$
\end{center}
and, in general:
\begin{equation}
\label{eq:symm}
\sum_{i=2}^{r} (-1)^{i} b_{ir}=
\sum_{j_1+2j_2+\cdots+rj_r=r} \frac{(j_1+\cdots +j_r)!}{j_1!\cdots j_r!}
(-b_1)^{j_1}\cdots(-b_r)^{j_r}.
\end{equation}

\subsection{Relation with the OS-algebra}
\label{subsec:bsas}

We now discuss the relation between the
Betti numbers of the resolution of $\k$ over $A$
and the numbers $a_j = \dim_{\k} \Tor_1^E(A,\k)_j$.

\begin{lem}
\label{lem:b2j}
We have: $b_{22}=\binom{b_1+1}{2}+a_2$, and
\begin{equation}
\label{eq:b2j}
b_{2j} = a_j, \quad\text{for  $j>2$}.
\end{equation}
\end{lem}
\begin{proof}
Clearly,
\begin{equation}
\label{eq:b22}
b_{22} = b_1^2-b_2 = b_1^2-\tbinom{b_1}{2} +a_2 = \tbinom{b_1+1}{2} + a_2.
\end{equation}

 From the change of rings spectral sequence associated to the composition
of ring maps $E \rightarrow A \twoheadrightarrow \k$  (see \S\ref{sec:ss}),
we obtain a $5$-term exact sequence
\begin{equation}
\label{eq:fiveterm}
   \Tor_2^{E}(\k,\k)\longrightarrow \Tor_2^{A}(\k,\k) \longrightarrow
\Tor_1^{E}(A,\k) \longrightarrow
\Tor_1^{E}(\k,\k) \longrightarrow \Tor_1^{A}(\k,\k) \longrightarrow 0.
\end{equation}
Now take free resolutions:
\begin{equation*}
\label{eq:triexact}
\begin{CD}
\cdots @>>> E^{\binom{b_1 +1}{2}}(-2) @>>> E^{b_1}(-1) @>>> E @>>>\k
\longrightarrow 0,
\\
\cdots @>>> \bigoplus_jA^{b_{2j}}(-j)@>>> A^{b_1}(-1) @>>> A @>>> \k
\longrightarrow 0,
\\
&&\cdots @>>> \bigoplus_jE^{a_{j}}(-j) @>>> E @>>> A\longrightarrow 0.
\end{CD}
\end{equation*}
The $\Tor$'s are graded, and we have the following ranks for
the respective graded pieces:
\begin{center}
$\begin{array}{c|ccccc}
\text{degree} & \Tor_2^{E}(\k,\k) & \Tor_2^{A}(\k,\k) &  \Tor_1^{E}(A,\k)
& \Tor_1^{E}(\k,\k)
& \Tor_1^{A}(\k,\k) \\
\hline
0 & 0 &0 &0 &0 &0\\
1 & 0 &0 &0 &b_1 &b_1\\
2 & \binom{b_1+1}{2} & b_{22} &a_2 &0 &0\\
3 & 0 & b_{23} &a_3 &0 &0\\
\vdots & \vdots & \vdots & \vdots &\vdots &\vdots\\
\ell & 0 & b_{2\ell} &a_{\ell} &0 &0
\end{array}$
\end{center}
The result follows at once.  (Notice that since
$b_{22} = \tbinom{b_1+1}{2} + a_2$, the sequence \eqref{eq:fiveterm} is also
exact on the left.)
\end{proof}

\subsection{Witt's formula}
\label{subsec:witt}

Recall that the lower central series of a finitely-generated group $G$ is
the chain of normal subgroups, $G=G_1 \ge G_2 \ge G_3\ge\cdots$,
with $G_k$ equal to the commutator $[G_{k-1},G]$, for $k>1$ (see \cite{MKS}
as a basic reference).  The successive
quotients of the series are finitely-generated abelian groups.
Let $\phi_k(G)=\rank G_k/G_{k+1}$ denote the rank
of the $k$-th LCS quotient.

The LCS ranks of a finitely-generated free group $F_n$ are given
by Witt's formula:
\begin{equation}
\label{eq:w1}
\phi_k(F_n)=\tfrac{1}{k}\sum_{d\, \mid k} \mu(d) n^{\frac{k}{d}},
\end{equation}
or, equivalently:
\begin{equation}
\label{eq:w2}
\prod_{k=1}^{\infty} (1-t^k)^{\phi_k} = 1-nt.
\end{equation}
For example:
\begin{equation}
\label{eq:w3}
\phi_1(F_n) = n,\quad
\phi_2(F_n) = \tbinom{n}{2},\quad
\phi_3(F_n) = 2\tbinom{n+1}{3},\quad
\phi_4(F_n) =  \tfrac{n^2(n^2-1)}{4}.
\end{equation}

\subsection{The first three $\phi_k$'s}
\label{sec:phi123}
We now compute the ranks of the first three LCS quotients
of an arrangement group, in terms of concrete combinatorial data.
The results of the previous section allow us to recover
a formula of Falk \cite{f}, but with a slightly new perspective.
This approach illustrates that high degree syzygies which appear
near the beginning of the resolution can provide information about low degree
syzygies ``further out''. Contrast this to the approach
of Peeva \cite{Pe}, where the Orlik-Solomon algebra is
truncated in degree three.
\begin{cor}
\label{cor:phi123}
The ranks of the first three LCS quotients of an arrangement group
are given by:
\begin{align}
\label{phi1}
\phi_{1}& = b_1 = \vert\rm{hyperplanes}\vert\\
\label{phi2}
\phi_{2}& =\binom{b_1}{2} - b_2 =  a_2\\
\label{phi3}
\phi_{3}& = b_3 -\binom{b_1}{3} +b_1\bigg(\binom{b_1}{2}-b_2\bigg) + a_3
  =  b_3 -\binom{b_1}{3} + b_1a_2+a_3.
\end{align}
\end{cor}
\begin{proof}
As explained in \S\ref{subsec:ext}, we have:
\begin{equation}
\label{eq:hilbwitt}
\prod_{k=1}^{\infty}(1-t^k)^{\phi_k}\cdot
\Hilb\Big(\bigoplus_i \Ext_A^i(\k,\k)_i,t\Big)=1.
\end{equation}
Expand both terms, using \eqref{eq:symm} and \eqref{eq:b2j}:
\begin{equation*}
\begin{split}
\prod_{k=1}^{\infty} (1-t^k)^{\phi_k}=
1-\phi_1 t +\big(\tbinom{\phi_1}{2}-\phi_2\big)t^2-
\big(\tbinom{\phi_1}{3}-\phi_1\phi_2+\phi_3\big)t^3+
\\
\big(\tbinom{\phi_1}{4}+\phi_1\phi_3-\tbinom{\phi_1}{2}\phi_2
+\tbinom{\phi_2}{2}-\phi_4\big)t^4+
\cdots
\end{split}
\end{equation*}
\begin{equation*}
\Hilb\Big(\bigoplus_i \Ext_A^i(\k,\k)_i,t\Big)=1+  b_1 t + (b_1^2-b_2)t^2 +
(b_1^3-2b_1b_2+b_3+a_3)t^3+\cdots
\end{equation*}
and solve.
\end{proof}

Formula \eqref{phi2} can be rewritten in terms of the M\"obius
function of the intersection lattice as
$\phi_2= \sum_{X\in L_2(\A)} \phi_2(F_{\mu(X)})$.
Formulas \eqref{phi3} and \eqref{eq:bp23}
show that $\phi_3$ equals the
number of linear syzygies on the degree two generators of $I$:
\begin{equation}
\label{eq:p3}
\phi_3=b'_{23}.
\end{equation}
This expression for $\phi_3$ (essentially the
same as the one obtained by Falk in \cite{f}), depends on more
subtle combinatorial data than just the M\"obius function of $L(\A)$.
Nevertheless, it was shown by Falk in \cite{Fa89} that,
for all $k\ge 3$,
\begin{equation}
\label{eq:phikbound}
\phi_k\ge \sum_{X\in L_2(\A)}\phi_k(F_{\mu(X)}).
\end{equation}
For $k=3$, this inequality also follows from  \eqref{eq:bpineq},
together with \eqref{eq:p3} and \eqref{eq:w3}.

\section{The change of rings spectral sequence}
\label{sec:ss}
In this section, we examine the change of rings spectral sequence
more deeply.  This spectral sequence is due to Cartan and
Eilenberg \cite{CE}; a good reference for the material presented
here is Eisenbud \cite{e}, appendix A.3
(see also McCleary \cite{McC}, p.~508).

Our goal is to compute the Betti numbers
$b_{ii}=\dim_{\k} \Tor_i^A(\k,\k)_i$, which will allow
us to solve for $\phi_k$. We concentrate on the case $k=4$,
since it illustrates well the general situation.
The Betti number $b_{44}$
can be obtained from $b_{34}$ and a Hilbert series computation as above,
or directly from the spectral sequence. We first tackle $b_{34}$,
since it corresponds to earlier terms in the spectral sequence.

\subsection{Change of rings}
\label{subsec:crss}
Suppose we have a composition of ring maps:
\begin{equation*}
E \rightarrow A \twoheadrightarrow \k.
\end{equation*}
By taking a free resolution $Q_{\bullet}$ for $\k$ over $E$ and a
free resolution $P_{\bullet}$ for $\k$ over $A$ and tensoring with $\k$,  we
obtain a first quadrant double
complex, which yields a spectral sequence
\begin{equation}
^2{E^{ij}}\cong \Tor_i^A\left(\Tor_j^E(A,\k),\k\right)
\Longrightarrow \Tor_{i+j}^E(\k,\k).
\end{equation}
Since $^1_{\hor}E^{i,j} = 0$ unless $i=0$, 
\begin{equation}
^2_{\hor}E^{i,j} =\: ^\infty_{\hor}E^{i,j}.
\end{equation}

In our situation (with $E\to A$ the canonical projection of the exterior
algebra onto the OS-algebra, and $A\twoheadrightarrow \k$ the projection
onto the residue field), we have
\begin{equation}
^2_{\hor}E^{0,j} = \Tor_j^E(\k,\k),
\end{equation}
and we know how to compute this term from the Koszul complex.
In particular, $\Tor_j^E(\k,\k)$
is nonzero only in degree $j$, and we have
\begin{equation}
\label{eq:torekk}
   \dim_{\k} \Tor_j^E(\k,\k)_j = \binom{b_1 + j -1}{j} .
\end{equation}
We already used this fact to show that the sequence 
\eqref{eq:fiveterm} of low degree terms was actually
left exact. To keep things easy to follow, we display the
$^2_{\ver}E^{i,j}$ terms:

{\small
\begin{equation*}
\xymatrixrowsep{16pt}
\xymatrixcolsep{12pt}
\xymatrix{
\Tor_2^E(A,\k) & \Tor_1^A(\Tor_2^E(A,\k),\k) &
\Tor_2^A(\Tor_2^E(A,\k),\k) & \Tor_3^A(\Tor_2^E(A,\k),\k)
\\
\\
\Tor_1^E(A,\k) &  \Tor_1^A(\Tor_1^E(A,\k),\k) &
\Tor_2^A(\Tor_1^E(A,\k),\k) \ar[uull]_(.35){d_2^{2,1}}
&  \Tor_3^A(\Tor_1^E(A,\k),\k)
\\
\\
\k & \Tor_1^A(\k,\k) &  \Tor_2^A(\k,\k)  \ar[uull]^(.66){d_2^{2,0}}
& \Tor_3^A(\k,\k)  \ar[uull]^(.66){d_2^{3,0}}
\ar@{..>}[uuuulll]_(.2){d_3^{3,0}}
}
\end{equation*}
}
\vspace*{6pt}

  From the diagram, we compute:
\begin{align*}
_{\ver}^{\infty}E^{2,0}&=\ker d_2^{2,0}
\\
_{\ver}^{\infty}E^{1,1}&=
\coker d_2^{3,0}
\\
_{\ver}^{\infty}E^{0,2}&=
(\coker d_2^{2,1})/d_3^{3,0}(\ker d_2^{3,0})
\end{align*}

Let $\tot$ be the total complex associated with the spectral sequence.
We then have
\begin{equation*}
\gr(H_k(\tot)) \cong \bigoplus\limits_{i+j=k}\null
_{\ver}^{\infty}E^{i,j},
\end{equation*}
and thus
\begin{equation}
\label{eq:tot}
\gr(H_k(\tot))_r = 0, \quad \text{for $r>k$}.
\end{equation}

The exactness of the sequence \eqref{eq:fiveterm} of terms of low degree
implies that both $_{\ver}^{\infty}E^{1,1}$ and $_{\ver}^{\infty}E^{0,2}$
vanish.   Hence:
\begin{align*}
\im d_2^{3,0}&=\Tor_1^A(\k,\Tor_1^E(A,\k)),
\\
\coker d_2^{2,1}&=d_3^{3,0}(\ker d_2^{3,0}).
\end{align*}

Tracing through the spectral sequence and using the equalities above,
we obtain exact sequences
\begin{equation*}
\label{eq:chase}
\begin{CD}
&0\\
&@VVV\\
&_{\ver}^{\infty}E^{3,0}\\
&@VVV\\
0\longrightarrow &\ker d_2^{3,0} @>>>\Tor_3^A(\k,\k)
@>{d_2^{3,0}}>>\Tor_1^A(\k,\Tor_1^E(A,\k))  \longrightarrow  0.\\
&@VV{d_3^{3,0}}V\\
&\coker d_2^{2,1}
\\
&@VVV\\
&0
\end{CD}
\end{equation*}

\vspace*{3pt}

Now recall that $\gr(H_3(\tot))$ vanishes in degree greater than
three, and so, in particular,
$\left(_{\ver}^{\infty}E^{3,0}\right)_4=0$.
We thus have proved the following theorem.
\begin{thm}
\label{thm:bound1}
The following sequence is exact:
\begin{equation}
\label{eq:exact}
0 \rightarrow (\coker d_2^{2,1})_4 \longrightarrow
\Tor_3^A(\k,\k)_4 \longrightarrow  \Tor_1^A(\k, \Tor_1^E(A,\k))_4
\rightarrow 0.
\end{equation}
\end{thm}
Since the dimension of $\Tor_1^A(\k, \Tor_1^E(A,\k))_4$ is $b_1 a_3$, this
contribution to $b_{34}=\dim_{\k} \Tor_3^A(\k,\k)_4$
is easy to understand.  We shall analyze the contribution
from $\dim_{\k} (\coker d_2^{2,1})_4$ next.

\subsection{Quadratic syzygies and a Koszulness test}
\label{subsec:ktest}
  We have already seen
that if $a_3 > 0$, then $A$ is not Koszul. The sequence \eqref{eq:exact}
gives another necessary condition for Koszulness.  Indeed, if
\begin{equation}
\label{eq:tor24}
\dim_{\k} \Tor_2^E(A,\k)_4 >
\dim_{\k}d_2^{2,1}(\Tor_2^A(\k,\Tor_1^E(A,\k)))_4,
\end{equation}
then $A$ is not a Koszul algebra.

In order to apply this Koszulness test, we need a way to determine
$(\im d_2^{2,1})_4$. This is the content of the following lemma.
First, recall some definitions.

Given elements
$f_1,\dots ,f_k\in E$, the {\em syzygy module} is a submodule
of a free module on $\epsilon_1,\dots,\epsilon_k$, consisting
of $E$-linear relations among the $f_i$'s. An (exterior) 
{\em Koszul syzygy} is a relation of the form 
$f_i\epsilon_j \pm  f_j\epsilon_i$. 
See \cite{e} for more details.

\begin{lem}
\label{lem:diff}
The image of $d_2^{2,1}$ in degree $4$ equals $\Kosz_4$,
the space of minimal quadratic
syzygies on the generators of $I_2$ which are Koszul.
\end{lem}
\begin{proof}

Take (minimal) free resolutions:
\begin{align*}
P_{\bullet}\colon \qquad
&0 \longleftarrow \k \longleftarrow A \longleftarrow
A^{b_1}(-1) \longleftarrow A^{\binom{b_1+1}{2}}(-2) \oplus A^{a_2}(-2)
\longleftarrow  \cdots
\\[3pt]
Q_{\bullet}\colon \qquad
&0 \longleftarrow \k \longleftarrow E \longleftarrow
E^{b_1}(-1) \longleftarrow E^{\binom{b_1+1}{2}}(-2)
\longleftarrow  \cdots
\end{align*}

\noindent
and form the double complex

\begin{equation*}
\xymatrix{
P_0\otimes(A\otimes Q_2)  &
P_1\otimes(A\otimes Q_2) \ar[d]^(.4){d} \ar[l]_{\delta} &
P_2\otimes(A\otimes Q_2) \\
P_0\otimes(A\otimes Q_1)  &
P_1\otimes(A\otimes Q_1)  &
P_2\otimes(A\otimes Q_1)
\ar[d]^{d} \ar[l]_{\delta} \ar@{-->}[ull]_(.4){d_2} \\
P_0\otimes(A\otimes Q_0)  &
P_1\otimes(A\otimes Q_0)  &
P_2\otimes(A\otimes Q_0)
}
\end{equation*}
\vspace*{4pt}

\noindent
The $d_2$ differential is just the differential from the snake
lemma. Since all the generators of $I_2$ may be
written as products of linear forms, the result follows from a
diagram chase.
\end{proof}

In view of this lemma, we may rephrase the non-Koszul
criterion~\eqref{eq:tor24}, as follows.

\begin{cor}
\label{cor:bound2}
If the minimal resolution of $A$ over $E$ has minimal quadratic
non-Koszul syzygies, then $A$ is not a Koszul algebra.
\end{cor}

The next lemma gives an upper bound on the
number of (minimal) quadratic Koszul syzygies,
solely in terms of the M\"{o}bius function of
the intersection lattice of the arrangement $\A$.

\begin{lem}
\label{lem:koszsyz}
Let $\Kosz_4$ be the the space of quadratic
Koszul syzygies on the generators of $I_2$.  Then:
\begin{equation}
\label{eq:ko}
\dim_{\k} \Kosz_4\le \sum _{(X, Y)\in \binom{L_2(\A)}{2}}
\binom{\mu(X)}{2} \binom{\mu(Y)}{2}.
\end{equation}
where $\tbinom{S}{2}$ denotes the set of unordered pairs of distinct
elements of a set $S$. Moreover, if $\A$ is MLS, then equality holds
in \eqref{eq:ko}.
\end{lem}
\begin{proof}

Quadratic Koszul syzygies between elements of $I_2$ which come from
the same dependent set are consequences of linear syzygies
(because the free resolution for pencils is linear).
Any two distinct elements of $L_2$, say $X$ and $Y$,
give rise to $\binom{\mu(X)}{2}  \binom{\mu(Y)}{2}$
quadratic Koszul syzygies, whence the bound.

In the MLS case, all linear syzygies come from pencils, hence any
of the above Koszul syzygies is {\em not} a consequence
of the linear syzygies.
\end{proof}

\subsection{Computing $b_{ii}$}
\label{subsec:bii}
The basic idea is to divide the problem up, using the exact sequence
\begin{equation}
   0 \longrightarrow \ker d_2^{i,0} \longrightarrow \Tor^A_i(\k,\k)_i
\longrightarrow \im d_2^{i,0} \longrightarrow 0.
\end{equation}

We start by analyzing the kernel of $d_2^{i,0}$,
illustrating the method in the $i=4$ case. Using the fact that the
differentials in the spectral sequence
are graded, we obtain the following diagram:
{\small
\begin{equation*}
\xymatrixrowsep{14pt}
\xymatrixcolsep{12pt}
\xymatrix{
**[ll]\Tor_3^E(A,\k)_4
\\
**[l]\bullet\qquad &**[r] \Tor_1^A(\k,\Tor_2^E(A,\k))_4
\\ \\
**[l]\bullet\qquad & **[l]\bullet &**[r] {\ \:}\Tor_2^A(\k,\Tor_1^E(A,\k))_4
\\ \\
**[l]\bullet\qquad & **[l]\bullet & **[l]\bullet & **[r]\Tor_4^A(\k,\k)_4
\ar[uul]_{d_2} \ar[uuuull]_(.63){d_3}  \ar[uuuuulll]^(.82){d_4}
}
\end{equation*}
}

\noindent
Now recall from \eqref{eq:tot} that $\gr(H_4(\tot))_r=0$ if $r>4$.
An analysis of the terms in $\bigoplus_{i+j=4} {_{\ver}^{\infty}E^{i,j}}$
shows that $_{\ver}^{\infty}E^{4,0}=\ker d_4$ is the only nonzero term and
\begin{equation*}
\dim_{\k}\, {_{\ver}^{\infty}E^{4,0}} =
\dim_{\k} \Tor_4^E(\k,\k)_4=\binom{b_1+3}{4}.
\end{equation*}

 From the diagram
\begin{equation*}
\label{eq:chase2}
\begin{CD}
&&0\\
&&@VVV\\
0 \longrightarrow\ker d_4 @>>>\ker d_3
@>>>\Tor_3^E(A,\k)_4 \longrightarrow 0\\
&&@VVV\\
&&\ker d_2\\
&&@VVV\\
&&\Tor_1^A(\k,\Tor_2^E(A,\k))_4 \\
&&@VVV\\
&&0
\end{CD}
\end{equation*}
we find that
\begin{equation*}
\dim_{\k}\ker d_2^{4,0}= \tbinom{b_1+3}{4} + b'_{34} +b_1b'_{23}.
\end{equation*}
In the general case, we obtain the following.
\begin{lem}
\label{lem:bii}
$\DS{
\dim_{\k} \ker d_2^{i,0}= \binom{b_1+i-1}{i}+\sum_{j=0}^{i-3}b_{jj}
b'_{i-j-1,i-j}.
}$
\end{lem}
Next, we analyze the image of $d_2^{i,0}$. From the diagram
\begin{equation*}
\label{eq:chase3}
\begin{CD}
&&&&0&\\
&&&&@VVV&\\
&&&&\big(\, ^3_{\ver}E^{i-2,1}\big)_i  \\
&&&&@VVV&\\
0\longrightarrow
\big(\im d_2^{i,0}\big)_i @>{d_2}>>\Tor_{i-2}^A(\k,\Tor_1^E(A,\k))_i
@>>>\big(\coker d_2^{i,0}\big)_i\ &\longrightarrow 0\\
&&&&@VVV&\\
&&&&\big(\im d_2^{i-2,1}\big)_i& \\
&&&&@VVV&\\
&&&&0&\\
\end{CD}
\end{equation*}
we obtain the following.
\begin{thm}
\label{thm:biidelta}
Let
$
\delta_{i}=\dim_{\k} d_2\big(\Tor_{i-2}^A(\k,\Tor_1^E(A,\k))\big)_i +
               \dim_{\k} \big(\, ^3_{\ver}E^{i-2,1}\big)_i  \,
$.
Then:
\begin{equation*}
b_{ii}+\delta_{i}= \binom{b_1+i-1}{i}+\sum_{j=0}^{i-2}b_{jj}
b'_{i-j-1,i-j}.
\end{equation*}
\end{thm}
As a consequence, we can compute the Betti numbers $b_{ii}$ inductively,
provided we know all the numbers $b'_{j,j+1}$ and $\delta_{j}$.
An easy check shows that
$\big(\, ^3_{\ver}E^{2,1}\big)_4=0$.
Thus, $\delta_4=\dim_{\k}\, \im \big( d_2^{2,1}\big)_4$.
 From Lemma~\ref{lem:diff}, we know this equals $\dim_{\k} \Kosz_4$.
We thus obtain the following.

\begin{cor}
\label{cor:b44}
Let $\delta_4$ be the number of quadratic Koszul  syzygies
on the generators of $I_2$.  Then:
\begin{equation}
\label{eq:b44}
b_{44}= \binom{b_1+3}{4}+\bigg(\binom{b_1+1}{2} + a_2\bigg)a_2+
b'_{34} +b_1b'_{23}-\delta_4
\end{equation}
\end{cor}

\subsection{Computing $\phi_4$}
\label{subsec:phi4}

 From formula \eqref{eq:syintro}, and the formulas from 
Corollary~\eqref{cor:phi123},
we find the following expression for the rank of the fourth LCS quotient:
\begin{equation}
\label{eq:p4}
\phi_4=\frac{b_1(b_1^3-6b_1^2-b_1-2)}{8}-\frac{b_1(3b_1+1)}{2}a_2-\binom{a_2+1}{2}
-b_1(a_3+b_3)+b_{44}.
\end{equation}
Combining \eqref{eq:b44} and \eqref{eq:p4}, we obtain the following.

\begin{thm}
\label{thm:phi4formula}
The rank of the fourth LCS quotient of an arrangement group $G$ is
given by:
\begin{equation}
\label{eq:phi4}
\phi_4(G)=\binom{a_2}{2}+b'_{34}-\delta_4.
\end{equation}
\end{thm}

Using the bounds for $b'_{34}$ and $\delta_4$ from
Lemmas \ref{lem:bprimebound} and \ref{lem:koszsyz},
together with formulas \eqref{eq:a2} and \eqref{eq:w3},
we recover Falk's lower bound \eqref{eq:phikbound}
for $\phi_k$, in the particular case $k=4$.

\section{Minimal linear strand arrangements}
\label{sect:nonres}

We now apply the machinery developed in the previous
sections to the class of MLS arrangements,
introduced in \S\ref{sec:resos}.  We start with
a connection with the resonance varieties.

\subsection{Local resonance}
\label{subsec:locres}
To an arbitrary arrangement $\A$, with Orlik-Solomon
algebra $A=H^*(M(\A),\k)$, Falk
associated in \cite{Fa97} a sequence of {\em resonance varieties},
defined as $\mathcal{R}^p(\A)=\{a\in E_1\mid H^p(A,\cdot a)\ne 0\}$.
Best understood is the variety $\mathcal{R}^1(\A)$, which
depends only on $L_{\le 2}(\A)$.  To each flat $X\in L_2(\A)$
there corresponds a ``local" component, $L_{X}$, of dimension
$\mu(X)$.  But, in general, there are other components in
$\mathcal{R}^{1}(\A)$, as illustrated by the braid arrangements.
The next result shows that non-trivial resonance cannot
happen for MLS arrangements.

\begin{thm}
\label{thm:locnonres}
If $\A$ is MLS, then all components of $\mathcal{R}^1(\A)$ are local.
\end{thm}

\begin{proof}
Suppose we have non-local resonance, i.e., there exist nonzero
$a, b\in E_1$ with $0\ne a\wedge b\in I_2$.  Say,
$I_2=(f_1,\dots ,f_k)$, so $a\wedge b=\sum c_i f_i$.
Then $0=a\wedge a \wedge b=a\wedge \sum c_i f_i=\sum c_ia\wedge f_i$.
Hence,
\begin{equation*}
\begin{pmatrix}c_1 a \\ \vdots\\ c_k a
\end{pmatrix}
\cdot \begin{pmatrix}f_1&\dots &f_k\end{pmatrix}=0,
\end{equation*}
and so we have a (linear) syzygy. Since the syzygy involves
$f_1=(e_r-e_s)\wedge (e_t-e_s)$ and, by assumption, all linear
syzygies are local,
the only variables that can appear in $a$ correspond
to hyperplanes incident on the flat $X\in L'_2(\A)$ containing $\{r,s,t\}$.
Now pick any $f_i$ corresponding to a flat $Y\in L'_2(\A)$ with $Y\ne X$.
Since $|X\cap Y|\le 1$, the element $a$ must be equal to some $e_u$,
which is impossible.
\end{proof}

\begin{remark}
Another, more indirect proof of this theorem can be assembled
from results of \cite{CSai}, \cite{CScv}.  Indeed, it can be
checked that the condition \eqref{a3cond} for an arrangement
to be MLS coincides with the condition given in \cite[Theorem 7.9]{CSai}
for the completion of the Alexander invariant to decompose as a direct sum of
``local" Alexander invariants. From this decomposition, one infers
that the variety $\mathcal{R}^1(\A)$ decomposes into local components,
see \cite[Theorem 3.9]{CScv}.
\end{remark}

\begin{remark}
The converse to Theorem~\ref{thm:locnonres} does not hold.
Indeed, if $\A$ is a complex realization of the MacLane
matroid $\mathtt{ML}_8$, then the Alexander invariant of $\A$ does not
decompose (and so $\A$ is not MLS), even though all the components of
$\mathcal{R}^1(\A)$ are local, see \cite[Example 8.6]{CSai} and
\cite[Example 10.7]{Su}.
\end{remark}

\subsection{LCS formula for MLS arrangements}

Recall that, for any arrangement, we have $\phi_1=b_1$,
$\phi_2=a_2$, and $\phi_3=b'_{23}$.  The assumption that
$\A$ is MLS means that $\phi_3$ attains the lower bound
\eqref{eq:phikbound} predicted by the local contributions, i.e.:
\begin{equation}
\label{a3c}
\phi_3=2\sum _{X\in L_2(\A)} \tbinom{\mu(X)+1}{3}=
\sum_{X\in L_2(\A)} \phi_3(F_{\mu(X)}).
\end{equation}
The next result shows that, under this assumption,
$\phi_4$ is given by an analogous formula.

\begin{thm}
\label{thm:addphi4}
If $\A$ is an MLS arrangement with group $G$, then:
\begin{equation}
\phi_4(G)=\sum_{X\in L_2(\A)} \phi_4(F_{\mu(X)}).
\end{equation}
\end{thm}

\begin{proof}
By Lemmas~\ref{lem:bprimebound} and \ref{lem:koszsyz}, we have:
\begin{align}
b'_{i,i+1}&= i\sum _{X\in L_2(\A)} \binom{\mu(X)+i-1}{i+1},
\\
\delta_4 &= \sum _{(X, Y)\in \binom{L_2(\A)}{2}}
\binom{\mu(X)}{2} \binom{\mu(Y)}{2}.
\end{align}
 From formula \eqref{eq:phi4}, we obtain:
\begin{equation}
\phi_4=\binom{a_2}{2}+3\sum _{X\in L_2(\A)} \binom{\mu(X)+2}{4}-
\sum _{(X, Y)\in \binom{L_2(\A)}{2}} \binom{\mu(X)}{2} \binom{\mu(Y)}{2}.
\end{equation}
Plugging in $a_2=\sum_{X\in L_2(\A)} \binom{\mu(X)}{2}$, we find:
\begin{equation}
\phi_4=\sum _{X\in L_2(\A)} \frac{\mu(X)^2 (\mu(X)^2-1)}{4}.
\end{equation}
Witt's formula ends the proof.
\end{proof}

Much of the above argument works for arbitrary $b_{ii}$ and $\phi_k$,
though the spectral sequence chase that went into proving
Corollary~\ref{cor:b44} becomes much more difficult.
We state the expected end result as a conjecture about the LCS ranks
of an MLS arrangement (the simplest case of the more general Resonance LCS
Conjecture from \cite{Su}), leaving the verification for the future.

\begin{conj}[MLS LCS formula]
\label{conj:localLCS}
If $G$ is the group of a minimal linear strand arrangement, then, for 
all $k\ge 2$:
\begin{equation}
\phi_k(G)=\sum_{X\in L_2(\A)} \phi_k(F_{\mu(X)}).
\end{equation}
In other words:
\begin{equation}
\prod_{k=1}^{\infty}(1-t^k)^{\phi_k}=
(1-t)^{b_1} \prod_{X\in L_2(\A)} \frac{1- \mu(X)\, t}{1-t} .
\end{equation}
\end{conj}

It follows from \cite[Prop.~3.12]{Fa89} that the conjecture
holds for product $3$-arrange\-ments.  Yet there are many
MLS arrangements which are not products, as illustrated by
the examples in \S\ref{subsec:ex}.

\subsection{A simple situation}
\label{subsec:simple}
We now look at what happens in a very simple situation,
namely, when $\A$ is MLS, and the M\"obius function on
$L_2(\A)$ has values only $1$ or $2$.
This occurs precisely when $\phi_3=2\phi_2$,
or, equivalently,
\begin{equation}
\label{eq:simple}
a_3+b_3=\binom{b_1}{3}- (b_1-2)a_2.
\end{equation}
In this situation, $b_{33}=\binom{b_1+2}{3}+(b_1+2)a_2$.
Moreover,
\begin{equation}
\label{eq:simpledelta}
b'_{i,i+1}=ia_2,\quad \text{and}\quad
\delta_4 = \binom{a_2}{2}.
\end{equation}
As a consequence, we find:
\begin{cor}
\label{cor:simplephi4}
If $\phi_3=2a_2$, then
$b_{44}= \tbinom{b_1+3}{4}+a_2\tbinom{b_1+3}{2}+\tbinom{a_2+1}{2}$, and so
\begin{equation*}
\phi_4=3a_2.
\end{equation*}
\end{cor}
Moreover, in this case Conjecture \ref{conj:localLCS}
takes the form:
$\phi_k(G)=a_2 \phi_k(F_2)$, or:
\begin{equation}
\label{eq:simpleLCS}
\prod_{k=1}^{\infty}(1-t^k)^{\phi_k}=(1-t)^{b_1-2a_2}  (1-2t)^{a_2}.
\end{equation}

\subsection{Examples}
\label{subsec:ex}

We conclude this section with a few examples, illustrating the
above formulas and conjectures.

\begin{example}[$\text{X}_3$ arrangement]
\label{ex:x3again}
The simplest example of an arrangement for which the
OS-algebra is not quadratic is the $\text{X}_3$ arrangement from
Example~\ref{ex:x3}.  We readily compute $\phi_3=2a_2=6$; thus,
the arrangement satisfies condition \eqref{eq:simple}.
In this case, formula \eqref{eq:simpleLCS}
predicts:  $\prod_{k=1}^{\infty}(1-t^k)^{\phi_k}=(1-2t)^{3}$.
\end{example}

\begin{example}[Fan arrangements]
\label{ex:fan}

Let $\A$ and $\A'$ be the pair of arrangements
considered by Fan in \cite{Fn} (see Figure~\ref{fig:fan}).
It is easy to see that both arrangements are MLS
(in fact, they satisfy condition \eqref{eq:simple}).

The arrangement $\A$ has group $G=F_1\times F_2\times F_2\times F_2$,
and thus, by Witt's formula, 
$\prod_{k=1}^{\infty}(1-t^k)^{\phi_k(G)}=(1-t)(1-2t)^3$.
The arrangement $\A'$ has group $G'=F_1\times G_0$,
where $G_0$ is the group of the $\text{X}_3$ arrangement.
Conjecture \ref{conj:localLCS} predicts
$\prod_{k=1}^{\infty}(1-t^k)^{\phi_k(G')}=(1-t)(1-2t)^3$.

In other words, the two arrangements should have the same
$\phi_k$'s (and thus, the same $b_{ii}$'s),  even though the respective
matroids are different.  Nevertheless,
the combinatorial difference is picked up by
the graded Betti numbers of the free resolutions of
the respective quadratic OS-algebras.
Indeed, $\dim_{\k} \Tor_2^{E} (\overline{A},\k)_4=3$, whereas
  $\dim_{\k} \Tor_2^{E} (\overline{A}',\k)_4=4$.
\end{example}

\begin{figure}
\subfigure{%
\begin{minipage}[t]{0.3\textwidth}
\setlength{\unitlength}{0.7cm}
\begin{picture}(5,3.2)(0.5,-0.5)
\put(3,3){\line(1,-1){3}}
\put(3,3){\line(-1,-1){3}}
\put(3,3){\line(0,-1){3}}
\multiput(0,0)(6,0){2}{\circle*{0.3}}
\multiput(1.5,1.5)(1.5,0){3}{\circle*{0.3}}
\multiput(3,3)(0,-3){2}{\circle*{0.3}}
\end{picture}
\end{minipage}
}
\setlength{\unitlength}{0.7cm}
\subfigure{%
\begin{minipage}[t]{0.3\textwidth}
\begin{picture}(5,3.2)(-1.0,-0.5)
\put(3,3){\line(-1,-1){3}}
\put(3,3){\line(0,-1){3}}
\put(1.5,1.5){\line(1,0){3}}
\multiput(0,0)(6,0){2}{\circle*{0.3}}
\multiput(1.5,1.5)(1.5,0){3}{\circle*{0.3}}
\multiput(3,3)(0,-3){2}{\circle*{0.3}}
\end{picture}
\end{minipage}
}
\caption{\textsf{The Fan matroids}}
\label{fig:fan}
\end{figure}

\begin{figure}
\setlength{\unitlength}{0.7cm}
\begin{picture}(5,3.8)(0.5,-0.5)
\put(3,3){\line(1,-1){3}}
\put(3,3){\line(-1,-1){3}}
\put(3,3){\line(0,-1){3}}
\put(1.5,1.5){\line(1,0){3}}
\put(0,0){\line(1,0){6}}
\multiput(0,0)(6,0){2}{\circle*{0.3}}
\multiput(1.5,1.5)(1.5,0){3}{\circle*{0.3}}
\multiput(3,3)(0,-3){2}{\circle*{0.3}}
\end{picture}
\caption{\textsf{The matroid of $\operatorname{X}_2$}}
\label{fig:x2}
\end{figure}

\begin{example}[Kohno arrangement]
\label{ex:x2}
An example of a quadratic OS-algebra which is not Koszul is given by
Kohno's $\operatorname{X}_2$ arrangement
(the corresponding matroid is depicted in Figure~\ref{fig:x2}).
It is readily seen that $P(t)=(1+t)(1+6t+10t^2)$ and $a_3=0$;
in particular, $A$ is quadratic.
Corollary~\ref{cor:phi123} gives
$\phi_3=6$; hence, the arrangement satisfies condition \eqref{eq:simple}.
Conjecture \ref{conj:localLCS} predicts:
$\prod_{k=1}^{\infty}(1-t^k)^{\phi_k}=\frac{(1-2t)^5}{(1-t)^3}$.
 From \eqref{eq:simpledelta}, we compute:
$b'_{24}=15$, and $\delta_4=10$.
Thus, $b_{34}=5$, and so the algebra $A$ is not Koszul.
\end{example}

\begin{example}[$9_3$ Configurations]
\label{ex:papp}
Let $\A$ and $\A'$ be realizations of the classical
$(9_3)_1$ and $(9_3)_2$ configurations.  The two arrangements have the
same M\"obius function, and thus, the same Poincar\'e polynomial:
$P(t)=(1+t)(1+8t+19t^2)$.  Even so, $\A'$ is MLS ($a_3=2$),
yet $\A$ is non-MLS ($a_3=4$).
\end{example}

\section{Graphic arrangements}
\label{sec:graphic}

There is a particularly nice interpretation of our
results and formulas in the case of graphic arrangements.

\subsection{Betti numbers and generators of OS-ideal}
\label{subsec:bettigraph}

Let $\G$ be a simple graph on $\ell$ vertices, and let $\A_{\G}$
be the corresponding graphic arrangement in $\C^{\ell}$, with complement
$M$.   The Orlik-Solomon algebra $A=H^*(M,\k)$ is a quotient of the
exterior algebra on generators $e_i$ corresponding to the edges of $\G$.
The Poincar\'e polynomial is given by
\begin{equation}
\label{eq:chromatic}
P(M,t)=(-t)^{\ell}\chi_{\G}(-t^{-1}),
\end{equation}
where $\chi_{\G}(t)$
is the {\em chromatic polynomial} of the graph, see \cite{ot}.
A direct algorithm for computing the Betti numbers $b_i$ (and explicit formulas
for the first four) may be found in Farrell \cite{Far}.
As for the number of minimal generators of the Orlik-Solomon ideal,
we have the following.

\begin{lem}
\label{lem:b2jgraph}
For all $j>2$,
\begin{equation}
a_j =|\text{chordless $(j+1)$-cycles}|.
\end{equation}
\end{lem}

\begin{proof}
We know the $m$-cycles are dependent sets, and give rise to degree $m-1$
elements of the Orlik-Solomon ideal. The correspondence is
one-to-one because a dependent set $e_1\cdots e_m$ gives a relation
\begin{equation*}
\prod_i (-1)^i e_1 \cdots \hat{e}_i \cdots e_m.
\end{equation*}
But the lead monomial of some lower degree element of $I$ divides
$e_2 \cdots e_m$ iff that monomial corresponds to a chord of
   $(e_1\cdots e_m)$, a contradiction.
\end{proof}

This result has been independently obtained by Cordovil and Forge \cite{CF}.
Notice that the lemma does not generalize directly to arbitrary arrangements.
For example, for a line configuration
with four lines through a point, there are four dependent triples. But one
of the OS relations can be written as a sum of the other three.
This situation does not occur with graphic arrangements.

\subsection{Koszulness and supersolvability}
\label{subsec:koss}
Stanley proved in \cite{St} that a graphic arrangement $\A_{\G}$
is free iff it is supersolvable iff the graph $\G$ is chordal (i.e.,
every circuit in $\G$ has a chord).
For a nice exposition we refer the reader to Edelman-Reiner \cite{ER}.
Using Lemma~\ref{lem:b2j} and Lemma~\ref{lem:b2jgraph}, we obtain:

\begin{thm}
\label{thm:graphkoszul}
A graphic arrangement is supersolvable if and only if
its OS-algebra is Koszul (in fact, quadratic).
\end{thm}

The forward implication holds for all arrangements
(cf.~\cite[Theorem 4.6]{SY}), but the converse is
not known in general (see \cite[Problem 6.7.1]{y}).
Previously, the converse was only known to hold
for hypersolvable arrangements (see \cite{JP}).

\subsection{The first three LCS ranks}
\label{subsec:phigraph}
Recall from \S\ref{subsec:bettigraph} that $a_i$ equals the number
of chordless $(i+1)$-cycles in $\G$, for $i>2$.   For each integer
$r\ge 1$, define
\begin{equation}
\label{eq:kappadef}
\kappa_{r-1}:= \abs{\{\G'\subseteq \G \mid \G'\cong K_r\}}
\end{equation}
to be the number of complete subgraphs on $r$ vertices,
so that $\kappa_0=\abs{\text{vertices}}=\ell$,
$\kappa_1=\abs{\text{edges}}=b_1$,
$\kappa_2=\abs{\text{triangles}}=a_2$, etc. In other words,
$(\kappa_0,\kappa_1,\dots, \kappa_{\ell-1})$ is the $f$-vector
of the ``clique complex" associated to $\G$.

Using this notation, we now answer Problem~1.4 from Falk's
recent survey \cite{Fa00}, which asks for a combinatorial
formula for $\phi_3$ in the graphic setting.

\begin{cor}
\label{cor:phigraphic}
For a graphic arrangement, the first three LCS ranks are given by:
\begin{center}
$\begin{array}{ccccc}
          \phi_{1}& = &b_1 &=& \kappa_1\\
          \phi_{2}& = &\binom{b_1}{2} - b_2 &=&  \kappa_2\\
          \phi_{3}& =& b_3 + b_1a_2+a_3-\binom{b_1}{3}& =
&2(\kappa_2+\kappa_3)
\end{array}$
\end{center}
\end{cor}

\begin{proof}
For a graphic arrangement, the first three Betti numbers are given by:
\begin{center}
$\begin{array}{ccc}
          b_{1} &= &\kappa_1\\
          b_{2}& =& \binom{b_1}{ 2} - a_2\\
          b_{3}& =& \binom{b_1 }{ 3} - b_1a_2 -a_3
            + 2(\kappa_2+\kappa_3)
\end{array}$
\end{center}
Now use Corollary~\ref{cor:phi123}.
\end{proof}

\subsection{Resolution of a graphic OS-algebra}
\label{subsec:graphares}

To proceed further, we need to analyze in more detail the
resolution of the Orlik-Solomon algebra over the exterior algebra,
in the case of a graphic arrangement. We start with an important example.

\begin{figure}[ht]
\setlength{\unitlength}{0.8cm}
\begin{picture}(5,4)(-0.5,-2)
\xygraph{!{0;<14.4mm,0mm>:<0mm,12mm>::}
[]*-{\bullet}
(
-^{\DS{e_0}}[ddr]*-{\bullet}(-^{\DS{e_1}}[ll]),
-_{\DS{e_2}}[ddl]*-{\bullet},
-^(.75){\DS{e_3}}[d]*-{\bullet}
(-_{\DS{e_4}}[dr],-^{\DS{e_5}}[dl])
)
}
\end{picture}
\caption{\textsf{The complete graph $K_4$}}
\label{fig:k4}
\end{figure}

\begin{example}
\label{ex:braidres}
The braid arrangement (in $\C^4$) is the graphic arrangement
associated to the complete graph $K_4$ (see Figure~\ref{fig:k4}).
The free resolution of the
Orlik-Solomon algebra $A$ as a module over the exterior algebra
$E=\bigwedge^*(e_0,\dots,e_5)$ begins:
\[
0 \longleftarrow A \longleftarrow E \xleftarrow{\:\partial_1\:}
E^{4}(-2)
\xleftarrow{\:\partial_2\:} E^{10}(-3) \longleftarrow  \cdots,
\]
where $\partial_1=
\begin{pmatrix}
\partial e_{145}&
\partial e_{235}&
\partial e_{034}&
\partial e_{012}
\end{pmatrix}$, and $\partial_2=$
{\Small
\[
\begin{pmatrix}
\xymatrixrowsep{2pt}
\xymatrixcolsep{2pt}
\xymatrix{
e_1-e_4 & e_1-e_5 & 0 & 0 & 0 & 0 & 0 & 0 & e_3-e_0 & e_2-e_0\\
0 & 0 & e_2-e_3 & e_2-e_5 & 0 & 0 & 0 & 0 & e_0-e_1 & e_0-e_4\\
0 & 0 & 0 & 0 & e_0-e_3 & e_0-e_4 & 0 & 0 & e_1-e_5 & e_2-e_5\\
0 & 0 & 0 & 0 & 0 & 0 & e_0-e_1 & e_0-e_2 & e_3-e_5 & e_4-e_5
}
\end{pmatrix}\!.
\]
}
\vspace*{1pt}

To see how this goes, write the generators of $I_2$ as
$f_i=\alpha_i \wedge \beta_i$ ($i=1,\dots,4$); for example,
$f_1=\partial e_{145}=(e_1-e_4)\wedge(e_1-e_5)$.
Each $f_i$ generates two ``local" linear syzygies,
$\alpha_i$ and $\beta_i$, which appear in columns
$2i-1$ and $2i$ of the matrix $\partial_2$.  Let
\begin{equation*}
\eta_1=e_0-e_1-e_3+e_5\quad  \text{and}\quad
\eta_2=e_0-e_2-e_4+e_5.
\end{equation*}
Note that $\eta_1\wedge \eta_2=f_1-f_2+f_3+f_4\in I_2$.
Thus, $\eta_1$ and $\eta_2$ belong to the resonance variety
$\mathcal{R}^1$. (In fact, $\eta_1$ and $\eta_2$ span an essential,
$2$-dimensional component of $\mathcal{R}^1$, corresponding to the
neighborly partition $\Pi=(05|13|24)$, see \cite{Fa97}.)
We now get two new linear syzygies, expressing the fact that
$\eta_i\wedge (\eta_1\wedge \eta_2)=0$ in $E_3$.  Reducing modulo
the ``local" syzygies, we obtain the last two columns of $\partial_2$.
It is easy to see that the columns of $\partial_2$ are linearly independent. 
From Corollary~\ref{cor:phigraphic}, we know that $b'_{23}=10$.  
Hence, the columns of $\partial_2$ form a complete set of linear 
syzygies on $I_2$.

Clearly, each of these $5$ pairs of linear syzygies generates
$i+1$ (independent) linear $i$-th syzygies.
For example, the pair $(\alpha_i,\beta_i)$ yields
$\alpha_i\wedge \alpha_i=\beta_i\wedge \beta_i=
\alpha_i\wedge\beta_i+\beta_i\wedge\alpha_i=0$
as linear second syzygies.
Hence, $b'_{i,i+1}=5i$, for all $i\ge 2$.

Finally, note that there are $\tbinom{4}{2}=6$ quadratic Koszul
syzygies on the generators of $I_2$, but that they are
all consequences of the linear syzygies on $I_2$.
For example:
\begin{equation*}
\begin{split}
\partial e_{034}\cdot\partial e_{145}-\partial e_{145}\cdot\partial e_{034}
=
(e_0-e_1)\cdot \big( (e_3-e_0)\partial e_{145}+(e_0-e_1)\partial e_{235}+ \\
(e_1-e_5)\partial e_{034}+(e_3-e_5)\partial e_{012} \big) +
(e_0-e_3)\cdot (e_1-e_4)\partial e_{145} +\\
(e_1-e_5)\cdot (e_0-e_4)\partial e_{034}+
(e_3-e_5)\cdot (e_0-e_1)\partial e_{012}.
\end{split}
\end{equation*}
It is easy to check that there are no other quadratic
syzygies on the generators of $I_2$, and thus $b'_{24}=0$.
\end{example}

Now let $\A=\A_{\G}$ be an arbitrary graphic arrangement.

\begin{lem}
\label{lem:gresbnd}
For $i\ge 2$,
\begin{equation*}
b_{i,i+1}' = i(\kappa_2+\kappa_3).
\end{equation*}
\end{lem}

\begin{proof}
 From Corollary~\ref{cor:phigraphic}, we know there are
$b'_{23}=\phi_3=2(\kappa_2+\kappa_3)$ linear first syzygies
on the generators of $I_2$. From the discussion in
Example~\ref{ex:braidres}, we know that
$2\kappa_2$ of those syzygies are local. Each non-local resonance
component associated to a $K_4$ subgraph generates an additional
pair of linear first syzygies, as in the example above.
We claim these $2\kappa_3$ linear syzygies are linearly independent.

Indeed, for there to be a linear dependence between a set
of such syzygies, some pair would have to overlap at a position.
Those two syzygies must involve a common element, $f$,
of $I_2$. By assumption, the syzygies come from distinct sub-$K_4$'s,
so these two $K_4$'s share a common triangle corresponding to
$f$. Now on each triangle of a $K_4$, the resonant syzygies do not
involve the variables of that triangle (see the above example). So
the two syzygies on $f$ involving distinct $K_4$'s involve different
sets of variables, hence are independent.  This contradiction
proves the claim.

Thus, since we know that the number of linear first syzygies is
$2(\kappa_2+\kappa_3)$, we have identified them all.
Each pair of local linear first syzygies associated to a triangle
generates $i+1$ linear $i$-th syzygies. This is also the case
for the linear syzygy pairs associated to the $K_4$ subgraphs, as
described in the example above. Such syzygies will be independent,
since they are supported only on the linear syzygy pair coming
from a $K_4$, so we are done.
\end{proof}

\begin{remark}
In general there can be interplay between the syzygies arising
from resonance components, and this will be reflected in the free resolution
(this happens for example in the case of the non-Fano arrangement).
But as noted above, such a phenomenon does not occur for graphic
arrangements.
\end{remark}

\begin{lem}
\label{lem:gdelbnd}
$ \delta_4\le \binom{\kappa_2}{2} -6(\kappa_3+\kappa_4)$.
\end{lem}

\begin{proof}
There are at most $\binom{\kappa_2}{2}$ quadratic Koszul syzygies.
As seen in Example~\ref{ex:braidres}, the linear syzygies of a
$K_4$-subgraph kill the $6$ Koszul syzygies corresponding to that subgraph.
A similar computation shows that the linear syzygies of a $K_5$-subgraph
kill $36$ Koszul syzygies.
Of those, $30$ are killed by the $5$ sub-$K_4$'s of the $K_5$,
and the remaining $6$ are killed by the $K_5$ itself, whence the bound.
\end{proof}

\subsection{The fourth LCS rank}
\label{subsec:graphphi4}
We now return to the computation of the LCS ranks
of an arrangement group.

\begin{thm}
\label{thm:gphi4}
For a graphic arrangement,
\begin{equation*}
\label{eq:gp4}
\phi_4 \ge 3\kappa_2+9\kappa_3+6\kappa_4.
\end{equation*}
Furthermore, if $\kappa_3=0$, then $\phi_{4} =3 \kappa_2$.
\end{thm}
\begin{proof}
The inequality follows from Lemmas~\ref{lem:gresbnd} and \ref{lem:gdelbnd},
together with Theorem~\ref{thm:phi4formula}. If $\kappa_3=0$, 
equality follows from
Corollaries~\ref{cor:simplephi4} and \ref{cor:phigraphic}.
\end{proof}

\begin{figure}
\setlength{\unitlength}{0.6cm}
\begin{picture}(5,4.8)(0,-2.5)
\xygraph{!{0;<12mm,0mm>:<0mm,12mm>::}
[]*-{\bullet}
(-[dd]*-{\bullet},-[dr]*-{\bullet}
(-[ur],-[dl]),-[rr]*-{\bullet}
(-[dd]*-{\bullet}
(-[ll],-[ul])))
}
\end{picture}
\caption{\textsf{A non-hypersolvable graph}}
\label{fig:nonhyper}
\end{figure}

\begin{example}
\label{ex:nonhyper}
Let $\G$ be the graph in Figure~\ref{fig:nonhyper}.
The arrangement $\A_{\G}$ is the simplest example
(in terms of number of vertices) of a graphic arrangement
which is not hypersolvable. In particular, no general
method could be applied to obtain a value for $\phi_4$.
On the other hand, Theorem~\ref{thm:gphi4} yields $\phi_4=12$.
\end{example}

\subsection{Graphic LCS formula}
\label{subsec:graphLCS}
Extensive computations suggest that equality holds in Lemma~\ref{lem:gdelbnd},
and hence in Theorem~\ref{thm:gphi4}. In fact, this seems to be part of a more
general pattern, which leads us to formulate the following conjecture.

\begin{conj}[Graphic LCS formula]
\label{conj:graphics}
If $\G$ is a graph with $\ell$ vertices, then:
\begin{equation}
\label{eq:GLCS}
\prod_{k=1}^{\infty} \left(1-t^k\right)^{\phi_k}=
\prod_{j=1}^{\ell-1} \left(1-jt\right)^{\DS{\sum_{s=j}^{\ell-1}(-1)^{s-j}
\tbinom{s}{j} \kappa_{s}
}}
\end{equation}
\end{conj}
Expanding both sides, \eqref{eq:GLCS} becomes equivalent
to the following sequence of equalities:
\begin{align*}
\phi_1&=\kappa_1\\
\phi_2&=\kappa_2\\
\phi_3&=2(\kappa_2+\kappa_3)\\
\phi_4&=3(\kappa_2+3\kappa_3+2\kappa_4)\\
\phi_5&=6(\kappa_2+5\kappa_3+8\kappa_4+4\kappa_5)\\
\phi_6&=9 \kappa_2 + 89 \kappa_3 + 260 \kappa_4 +
300 \kappa_5 + 120 \kappa_6\\
\vdots&\qquad\vdots\\
\phi_k&=\sum_{j=1}^{k}\sum_{s=j}^{k} (-1)^{s-j} \binom{s}{j}
\kappa_{s} \phi_k(F_j)
\end{align*}
By Corollary~\ref{cor:phigraphic}, these equalities hold
up to $k=3$.  By Theorem~\ref{thm:gphi4},
equality also holds for $k=4$, provided $\kappa_3=0$
(otherwise, we only know there is an inequality in one direction).

\begin{remark}
\label{rem:chen}
For graphic arrangements, the resonance formula for the ranks
of the Chen groups, $\theta_k(G)=\rank\, \gr_k(G/G'')$
(conjectured in \cite{Su}) gives:
\begin{equation}
\theta_k=(k-1)(\kappa_2+\kappa_3),\quad \text{for $k\ge 3$.}
\end{equation}
The resonance LCS formula (also conjectured in \cite{Su}) is based upon
the assumption that $\phi_4=\theta_4$.  In the case of graphic arrangements,
formula \eqref{eq:GLCS} would imply:
\begin{equation}
\phi_3=2\phi_2  \Longleftrightarrow \phi_4=\theta_4  \Longleftrightarrow
\kappa_3=0.
\end{equation}
Hence, if $\G$ contains no $K_4$ subgraphs, then all the above conjectures
reduce to:
\begin{equation}
\theta_k=\kappa_2 \theta_k(F_2)\ \text{ and }\
\phi_k=\kappa_2 \phi_k(F_2),
\quad \text{for all $k\ge 2$}.
\end{equation}
\end{remark}

\begin{remark}
\label{rem:chordal}
The Graphic LCS conjecture is true for chordal graphs. Indeed,
if $\G$ is a chordal graph, then its chromatic polynomial is given by
\begin{equation}
\label{eq:chordchi}
\chi_{\G}(t)=t^{\kappa_0} \prod_{j=1}^{\kappa_0-1}
\left(1-jt^{-1}\right)^{\sum_{s=j}^{\kappa_0-1}(-1)^{s-j}
\tbinom{s}{j} \kappa_{s}},
\end{equation}
see \cite{BKWW}.  Since in this case $\A_{\G}$ is supersolvable,
formula~\eqref{eq:GLCS} follows at once from \eqref{eq:chromatic}
and the Falk-Randell LCS formula~\eqref{lcsformula}.
We warmly thank Yuri Volvovski for bringing  formula~\eqref{eq:chordchi}
to our attention.
\end{remark}

\bibliographystyle{amsalpha}

\end{document}